\documentclass[preprint,11pt]{elsarticle}
\usepackage{amssymb,amsmath,amsthm,amsfonts,amscd}
\usepackage{graphicx,color}
\usepackage[normalem]{ulem} 
\usepackage{soul} 
\usepackage{cancel} 
\usepackage{hyperref}
\definecolor{ForestGreen}{rgb}{0.15,0.416,0.18}
\definecolor{EgyptBlue}{rgb}{0.063,0.2,0.65}
\hypersetup{
	colorlinks=true,
	linkcolor=EgyptBlue,
   citecolor=EgyptBlue,
   urlcolor=ForestGreen
}
\journal{xxxxx}
\newtheorem{theorem}{Theorem}[section]

\newtheorem{definition}{Definition}[section]
\newtheorem{remark}{Remark}[section]

\newcommand{\cqdf}{\hfill\rule{5pt}{5pt}}

\newcommand{\disp}{\displaystyle}

\def\Om{\Omega}

\def\re{\mathbb{R}}

\def\jnt{\displaystyle\int}
\def\jjnt{\jnt\!\!\!\jnt}

\def\beq{\begin{eqnarray}}
\def\eeq{\end{eqnarray}}
\def\beqa{\begin{eqnarray*}}
\def\eeqa{\end{eqnarray*}}
\def\beqn{\begin{equation}}
\def\eeqn{\end{equation}}
\def\bar{\begin{array}}
\def\ear{\end{array}}
\begin{document}

\begin{frontmatter}

\title{Approximate controllability for a one-dimensional wave equation with the fixed endpoint control}
\author[IJesus]{Isa\'{i}as Pereira de Jesus\corref{cor1}}
\ead{isaias@ufpi.edu.br}

\address[IJesus]{ Departamento de Matem\'{a}tica, Universidade Federal do Piau\'{i}, 64049-550, Teresina, PI, Brazil}
\cortext[cor1]{Corresponding author. Phone: +55 (86) 3215-5835}

\begin{abstract}
This paper is devoted to the study of the  approximate controllability for a one-dimensional wave equation in domains with moving boundary. This equation models
the motion of a string where an endpoint is fixed and the other one is moving.  When the
speed of the moving endpoint is less than the characteristic speed, the
controllability of this equation is established. We present the following results: the existence
and uniqueness of Nash equilibrium, the approximate controllability  with respect
to the leader control, and the optimality system for the leader control.
\end{abstract}


\end{frontmatter}

\section{Introduction}

The development of science and technology has motivated many branches of control theory. Initially, in the classical
control theory, we encountered problems where a system must reach a predetermined target by the action of a single control,
for example, find a control of minimum norm such that the design specifications are met. To the extent that more realistic situations
were considered, it was necessary to include several different(and even contradictory) control objectives, as well as develop
theory that would handle the concepts of multi-criteria optimization, where optimal decisions need to be taken in the presence of trade-offs
between these different objectives. There are many points of view to deal with multi-objective problems. Notions of economics and game
theory were introduced in the works of J. Nash \cite{N}, V. Pareto \cite{P} and H. von Stackelberg \cite{S}, where
each has a particular philosophy to solve these problems.

According to the formulation introduced by H.von Stackelberg \cite{S}, we assume the presence of various local controls, called {\it followers}
which have their own objectives, and a global control, called {\it leader}, with a different goal from the rest of the players (in the case, the followers).
The general idea of this strategy is a game of hierarchical nature, where players compete against each other, so that the leader makes the first move
and then followers react optimally to the action of the leader. Since many followers are present and each has a specific objective, it is intended that
these are in Nash equilibrium.

This paper was inspired by ideas of J.-L. Lions
\cite{L1}, where we  investigate a  similar question of hierarchical control
employing the Stackelberg strategy in the case of time dependent domains.

Up to date, in the context of partial differential equations(PDEs), there are several papers related to this topic.
The papers by Lions  \cite {L10, L12},  the author gives some results concerning
Pareto and Stackelberg strategies, respectively. The paper by  D\'iaz and Lions
\cite {DL}, the approximate controllability of a system is established following a
Stackelberg-Nash strategy and the extension in D\'iaz \cite{D1}, that provides a
characterization of the solution by means of Fenchel-Rochafellar duality theory.
In \cite {RA1, RA2},  Glowinski, Ramos and Periaux   analyze the Nash
equilibrium for constraints given by linear parabolic and Burger's equations from
the mathematical and numerical viewpoints. The Stackelberg-Nash strategy for
the Stokes systems has been studied by Gonz\'alez, Lopes and Medar in \cite
{GO}. In Limaco, Clark and Medeiros \cite {LI}, the authors present the
Stackelberg-Nash equilibrium in the context of linear heat equation in non-cylindrical domains. The paper by Araruna, Fern\'andez-Cara and Santos
\cite{AR} ,  the authors developed the first hierarchical results within the exact controllability framework for a parabolic equation.
Finally, we can mention the paper by Ramos and Roubicek \cite{RR},  where the existence of a Nash equilibrium is
proved for a nonlinear distributed parameter predator-prey system and a
conceptual approximation algorithm is proposed and analyzed.

In this paper, we present the following results: the existence and uniqueness of
Nash equilibrium, the approximate controllability  with respect to the leader
control, and the optimality system for the leader control.

\subsection{Organization of the paper}

The remainder of the paper is organized as follows. In Section \ref{sec2}, we present the
problem. Section \ref{sec3} is devoted to establish the existence and uniqueness of Nash equilibrium. In Section \ref{sec4}, we study the approximate controllability with respect to the leader control. Finally, in the Section \ref{sec5} we present
the optimality system for the leader control.

\section{Problem formulation}\label{sec2}
\setcounter{equation}{0}

As in \cite{Cui10}, given $ T > 0 $, we consider the non-cylindrical domain defined
by
$$\displaystyle \widehat{Q}= \left\{ (x,t) \in \re^2;\; 0 < x <  \alpha _k(t),\; \; t \in (0,T) \right\},$$
where
$$ \alpha_k(t) = 1 + kt, \;\;\;\;\;\; \; 0 < k < 1.$$

Its lateral boundary  is defined by $\displaystyle \widehat{\Sigma} =\widehat{\Sigma}_0
\cup \widehat{\Sigma}_0^*$, with
$$\displaystyle \widehat{\Sigma}_{0} = \{(0, t);\; t \in (0,T) \} \;\;\;\;\mbox{ and } \;\;\;\; \widehat{\Sigma}_{0}^* =\widehat{\Sigma} \backslash \widehat{\Sigma}_0= \{(\alpha_k(t),t);\;  t \in (0,T) \}.$$
We also represent by $\Om_t$ and $\Om_0$ the intervals $\displaystyle (0, \alpha_k(t))$
and $\displaystyle (0, 1)$, respectively.

Consider the following  wave equation in the non-cylindrical domain
$\widehat{Q}$:
\begin{equation} \label{eq1.3}
\left|
\begin{array}{l}
\displaystyle u'' - u_{xx} = 0 \ \ \mbox{ in } \ \ \widehat{Q},\\[11pt]
\displaystyle u(x,t) = \left\{
\begin{array}{l}
\widetilde{w}(t) \ \ \mbox{ on } \ \ \widehat{\Sigma}_{0},\\[11pt]
0 \ \ \mbox{ on } \ \ \widehat{\Sigma}_0^*,
\end{array}
\right.\\[11pt]
\displaystyle u(x,0) = u_0(x), \;\; u'(x,0) = u_1(x) \ \mbox{ in }\; \Om_0,
\end{array}
\right.
\end{equation}
where $u$ is the state variable, $\widetilde{w}$ is the control variable and \linebreak $( u_0(x),
u_1(x)) \in L^2(0,1) \times H^{-1}(0,1)$. By $\displaystyle u'=u'(x,t)$ we represent the
derivative $\displaystyle \frac{\partial u}{\partial t}$ and by  $\displaystyle u_{xx}=u_{xx}(x,t)$ the second order partial derivative $\displaystyle \frac{\partial^2 u}{\partial x^2}$. Equation \eqref{eq1.3} models the motion of a string where an endpoint is fixed and the other one is moving. The constant $k$ is called the speed of the moving endpoint.

In \cite{Cui10}, instead of transforming the problem \eqref{eq1.3} from a non-cylindrical domain into a cylindrical domain, Cui et al. studied the controllability problem directly in a non-cylindrical domain, when the control is put on the fixed endpoint.

The control system of this paper is the same as that of \cite {Je} and \cite {Je1}. But motivated by \cite {Cui10}, we extend the result in \cite{Je} and \cite{Je1}, and the controllability result is obtained when $k\in (0,1).$

It would be quite interesting to study the controllability for   \eqref{eq1.3} in the
non-cylindrical domain  $\widehat{Q}$ for the case of $k > 1$. However, it seems
very difficult and remains to be done.

In practical situations, many processes evolve in domains whose boundary has moving parts. A simple model is, e.g., the interface of an ice-water mixture when the temperature rises. To study the controllability problem of wave equations with moving boundary or free boundary is very significant. In addition, for $0 < k < 1 $, \eqref{eq1.3} admits a unique solution in the sense of transposition. But for $k > 1,$ \eqref{eq1.3} admits a unique solution only when more boundary conditions are described;  see Cui et al. \cite{Cui10}.

In spite of a vast literature on the controllability problems of the wave equation in cylindrical domains, there are only
a few works dealing with non-cylindrical case. We refer to  \cite {Cui10, Je, Je1, Cui, Ar, Bar, Cui1, Je10, Mi1, Mi} for some known results in this direction.

In this article, motivated by the arguments contained in the work of J.-L. Lions
\cite{L1}, we investigate a  similar question of hierarchical control for equation
\eqref{eq1.3}, employing the Stackelberg strategy in the case of time
 dependent domains.

Following the work of J.-L. Lions \cite{L1}, we divide $\displaystyle \widehat{\Sigma}_0$
into two disjoint  parts
\begin{equation}\label{decomp0}
\displaystyle \widehat{\Sigma}_0=\displaystyle \widehat{\Sigma}_1 \cup \displaystyle \widehat{\Sigma}_2,
\end{equation}
 and consider
\begin{equation} \label{decomp}
\displaystyle \widetilde{w}=\{\widetilde{w}_1, \widetilde{w}_2\}, \;\; \widetilde{w}_i=\mbox{control function in } \; L^2(\widehat{\Sigma}_i), \;i=1,2.
\end{equation}

We can also write
\begin{equation} \label{decomp 2.A}
\displaystyle \widetilde{w}= \widetilde{w} _1 + \widetilde {w}_2, \; \mbox{ with } \; \displaystyle {\widehat \Sigma}_0=\displaystyle {\widehat \Sigma}_1 = \displaystyle {\widehat \Sigma}_2.
\end{equation}

Thus, we observe  that the system \eqref{eq1.3} can be rewritten as follows:
\begin{equation} \label{eq1.3.1}
\left|
\begin{array}{l}
\displaystyle u''- u_{xx} = 0 \ \ \mbox{ in } \ \ \widehat{Q},\\[5pt]
\displaystyle u(x,t) = \left\{
\begin{array}{l}
\widetilde{w}_1(t) \ \ \mbox{ on } \ \ \widehat{\Sigma}_{1},\\
\widetilde{w}_2(t) \ \ \mbox{ on } \ \ \widehat{\Sigma}_{2},\\
0 \ \ \mbox{ on } \ \ \widehat{\Sigma}\backslash \widehat{\Sigma}_0,
\end{array}
\right.\\[5pt]
\displaystyle u(x,0) = u_0(x), \;\; u'(x,0) = u_1(x) \ \mbox{ in }\; \Om_0.
\end{array}
\right.
\end{equation}
In the decomposition \eqref{decomp0}, \eqref{decomp} we establish a hierarchy.
We think of $\widetilde{w}_1$ as being the ``main'' control, the leader,  and we
think of $\widetilde{w}_2$ as the follower,  in Stackelberg terminology.

We assume that the objective of the leader is of controllability type. In fact, the leader wants to drive the state $u$ at the final time $T$ as close as possible to a wished state $u^0$, without a big cost for the control $\widetilde{w}_1$ . On the other hand, the main objective of the follower is to hold the state $u$ near to a desired state $\widetilde{u}_2(x,t)$ for all time $t \in (0,T)$, without a big cost for the control $\widetilde{w}_2$. In order to combine these two different objectives, the Stackelberg strategy is used (cf. J.-L. Lions \cite{L1}).

Associated with the solution $u=u(x,t)$ of \eqref{eq1.3.1}, we will consider the
(secondary) functional
\begin{equation}\label{sfn}
\displaystyle \widetilde{J}_2(\widetilde{w}_1, \widetilde{w}_2) =
\displaystyle\frac{1}{2} \displaystyle \jjnt_{\widehat{Q}} \left(u(\widetilde{w}_1, \widetilde{w}_2)-\widetilde{u}_2\right)^2 dx dt +
\displaystyle\frac{\widetilde{\sigma}}{2} \int_{\widehat{\Sigma}_2} \widetilde{w}_2^2\;d\widehat{\Sigma},
\end{equation}
and the (main) functional
\begin{equation}\label{mfn}
\displaystyle \widetilde{J}(\widetilde{w}_1) =\frac{1}{2}\int_{\widehat{\Sigma}_1} \widetilde{w}_1^2\;d\widehat{\Sigma},
\end{equation}
where $\widetilde{\sigma}>0$ is a constant and $\widetilde{u}_2$ is a given
function in $L^2(\widehat{Q}).$

\begin{remark}\label{bdfnc} From the regularity and uniqueness of the solution  of \eqref{eq1.3.1} (see Remark \eqref{rsol}) the cost functionals $\displaystyle \widetilde{J}_2$ and $\displaystyle \widetilde{J}$ are well defined.
\end{remark}
The control problem that we will consider is as follows: the follower $\displaystyle
\widetilde{w}_2$ assumes that  the leader $\displaystyle \widetilde{w}_1$ has made  a
choice. Then, it tries  to find an equilibrium of the cost  $\widetilde{J}_{2}$ , that
is, it looks for  a control $\displaystyle \widetilde{w}_2=\mathfrak{F}(\widetilde{w}_1)$
(depending on $\displaystyle \widetilde{w}_1$), satisfying:
\begin{equation}\label{n1.1}
\displaystyle\widetilde{J}_2(\widetilde{w}_1, \widetilde{w}_2) \leq \widetilde{J}_2(\widetilde{w}_1,  \widehat{w}_2), \;\;\; \forall\; \widehat{w}_2  \in L^2(\widehat{\Sigma}_2).
\end{equation}

The control $ \widetilde{w}_2 $, solution of \eqref{n1.1}, is called Nash
equilibrium for the cost $\widetilde{J}_2 $ and it depends on $ \widetilde{w}_1$
(cf. Aubin \cite{A}).

\begin{remark}\label{r2} In another way, if the leader $ \widetilde{w}_1$ makes a choice, then the follower $ \widetilde{w}_2$ makes also a choice, depending on  $\widetilde{w}_1$, which minimizes the cost
 $\widetilde{J}_2$, that is,
\begin{equation}\label{son}
\displaystyle \widetilde{J}_2(\widetilde{w}_1, \widetilde{w}_2)= \inf_{\widehat{w}_2 \in L^2(\widehat{\Sigma}_2)} \widetilde{J}_2(\widetilde{w}_1, \widehat{w}_2).
\end{equation}

This is equivalent to \eqref{n1.1}. This process is called Stackelberg-Nash
strategy; see \textup{D\'iaz and Lions \cite{DL}.}
\end{remark}

After this, we consider the state $\displaystyle u\left(\widetilde{w}_1,
\mathfrak{F}(\widetilde{w}_1)\right)$ given by the solution of
\begin{equation} \label{eq1.3.1F}
\left|
\begin{array}{l}
\displaystyle u''- u_{xx} = 0 \ \ \mbox{ in } \ \ \widehat{Q},\\[3pt]
\displaystyle u(x,t) = \left\{
\begin{array}{l}
\widetilde{w}_{1} \ \ \mbox{ on } \ \ \widehat{\Sigma}_{1},\\
\mathfrak{F}(\widetilde{w}_1) \ \ \mbox{ on } \ \ \widehat{\Sigma}_{2},\\
0 \ \ \mbox{ on } \ \ \widehat{\Sigma}\backslash \widehat{\Sigma}_0,
\end{array}
\right.\\[3pt]
\displaystyle u(x,0) = u_0(x), \;\; u'(x,0) = u_1(x) \ \mbox{ in }\; \Om_0.
\end{array}
\right.
\end{equation}

We will look for any optimal control $\displaystyle \widetilde{w_1}$ such that
\begin{equation}\label{ocn}
 \displaystyle \widetilde{J}(\widetilde{w}_{1}, \mathfrak{F}(\widetilde{w}_1))= \inf_{\overline{w}_1 \in L^2(\widehat{\Sigma}_1)} \widetilde{J}(\overline{w}_1, \mathfrak{F}(\widetilde{w}_1)),
\end{equation}
subject to the following restriction of the approximate controllability type
\begin{equation}\label{apcn}
\displaystyle \left(u(x, T; \widetilde{w}_{1}, \mathfrak{F}(\widetilde{w}_1)), u'(x, T; \widetilde{w}_{1}, \mathfrak{F}(\widetilde{w}_{1}))\right)\in B_{L^2(\Omega_t)}(u^0,\rho_0) \times B_{H^{-1}(\Omega_t)}(u^1,\rho_1),
\end{equation}
where $\displaystyle B_X(C, r)$ denotes the ball in $X$ with centre $C$ and radius $r$.

To explain this optimal problem, we are going to consider the following
sub-problems:

\textbf{$\bullet$ Problem 1} Fixed any leader control $\widetilde{w}_1$,  find the
follower control $\displaystyle \widetilde{w}_2=\mathfrak{F}(\widetilde{w}_1)$ (depending
on $\displaystyle \widetilde{w}_1$) and the associated state $u$, solution of
\eqref{eq1.3.1} satisfying the condition \eqref{son} (Nash equilibrium)  related to
$\widetilde{J}_2$, defined in \eqref{sfn}.

\textbf{$\bullet$ Problem 2} Assuming that the existence of the Nash equilibrium
$\displaystyle \widetilde{w}_2$ was proved, then when $\widetilde{w}_{1}$ varies in
$L^2(\widehat{\Sigma}_1)$, prove that the solutions $\displaystyle \left(u(x, t;
\widetilde{w}_1, \widetilde{w}_2), u'(x, t; \widetilde{w}_1, \widetilde{w}_2)\right)$ of
the state equation \eqref{eq1.3.1}, evaluated at $t = T$, that is, $\displaystyle \left(u(x, T;
\widetilde{w}_1, \widetilde{w}_2), u'(x, T; \widetilde{w}_1, \widetilde{w}_2)\right)$,
generate a dense subset of \break $L^2(\Omega_t) \times H^{-1}(\Omega_t)$.

\begin{remark}\label{r1} By the linearity of system \eqref{eq1.3.1F}, without loss of generality we may assume that $u_0=0=u_1$.
\end{remark}

\begin{remark}\label{rsol} Using similar technique as  in \textup{\cite{Mi}}, we can prove the  following:
For each $u_0 \in L^2(0,1)$, $u_1 \in H^{-1}(0,1)$ and $\displaystyle \widetilde {w}_i \in
L^2(\widehat{\Sigma}_i), \;i=1,2,$ there exists exactly one  solution $\displaystyle u$ to
\eqref{eq1.3.1} in the sense of a transposition, with  $\displaystyle u \in C\big( [0,T]
;L^2(\Om_t) \cap C^{1}\big( [0,T] ; H^{-1}(\Om_t))$. Thus, in particular, the cost
functionals $\displaystyle \widetilde{J}_2$ and $\displaystyle \widetilde {J}$ are well defined.
\end{remark}
\section{Nash equilibrium}\label{sec3}
\setcounter{equation}{0}
In this section, fixed any leader control $\displaystyle \ w_1 \in L^2(\widehat{\Sigma}_1)$ we
determine the existence and uniqueness of solutions to the problem
\begin{equation} \label{eq3.10}
\begin{array}{l}
\displaystyle\inf_{\widetilde{w}_2 \in L^2(\widehat{\Sigma}_2)}J_2(\widetilde{w}_1,\widetilde{w}_2),
\end{array}
\end{equation}
and a characterization of this solution in terms of an adjoint system.

In fact, this is a classical type problem in the control of distributed systems (cf.
J. -L. Lions \cite{L3}). It admits an unique solution
\begin{equation} \label{eq3.11}
\displaystyle \widetilde{w}_2 = \mathfrak{F}(\widetilde{w}_1).
\end{equation}
The Euler - Lagrange equation for problem  \eqref{eq3.10} is given by
\begin{equation} \label{eq3.21}
\int_{0}^{T}\int_{\Omega_t}(u - \widetilde{u}_2)\widehat{u}dxdt + \widetilde{\sigma}\int_{\widehat \Sigma_2}\widetilde{w}_2\widehat{w}_2d\ {\widehat \Sigma} = 0, \;\;\forall\, \widehat{w}_2 \in L^2(\widehat \Sigma_2),
\end{equation}
where $\widehat{u}$ is solution of the following system
\begin{equation} \label{eq3.22}
\left|
\begin{array}{l}
\displaystyle \widehat u'' - u_{xx} = 0 \ \ \mbox{ in } \ \ \widehat{Q},\\ [7pt]
\displaystyle \widehat{u} = \left\{
\begin{array}{l}
0 \ \ \mbox{ on } \ \ \widehat{\Sigma}_1, \\
\widehat{w}_2 \ \ \mbox{ on } \ \ \widehat{\Sigma}_2, \\
0 \ \ \mbox{ on } \ \ \widehat{\Sigma} \backslash \left({\widehat\Sigma_1} \cup {\widehat \Sigma_2}\right),
\end{array}
\right. \\[13pt]
\displaystyle \widehat{u}(x,0) = 0, \; \widehat{u'}(x,0) = 0, \;\; x \in \Om_t.
\end{array}
\right.
\end{equation}
In order to express \eqref{eq3.21} in a convenient form, we introduce the adjoint
state defined by
\begin{equation}\label{sac}
\left|
\begin{array}{l}
\displaystyle p'' - p_{xx} = u - \widetilde{u}_2  \ \mbox{ in } \ \ \widehat{Q}, \\[5pt]\displaystyle
p(T) = p'(T) = 0, \;\; x \in \Om_t, \\[5pt]\displaystyle
p = 0 \ \mbox{ on } \ \widehat \Sigma,
\end{array}
\right.
\end{equation}
Multiplying \eqref{sac} by $\widehat{u}$ and integrating by parts, we find
\begin{equation} \label{eq3.33}
\int_{0}^{T}\int_{\Om_t} (u - \widetilde{u}_2)\widehat{u}\,dx\,dt + \int_{\widehat \Sigma_2} p_x\,\widehat{w}_2\,d{\widehat \Sigma} = 0,
\end{equation}
so that \eqref{eq3.21} becomes
\begin{equation}\label{ci}
\displaystyle p_x= \widetilde{\sigma} \,\widetilde{w}_2 \ \ \mbox{ on } \ \ \widehat \Sigma_2.
\end{equation}
We summarize these results in the following theorem.
\begin{theorem}\label{teN} For each $\displaystyle \widetilde {w}_1 \in L^2(\Sigma_1)$ there exists a unique Nash equilibrium $\displaystyle \widetilde{w}_2$ in the sense of \eqref{son}. Moreover, the follower $\displaystyle \widetilde{w}_2$ is given by
\begin{equation}\label{cseg}
\displaystyle \displaystyle \widetilde{w}_2 = \mathfrak{F}(\widetilde {w}_1)=\frac{1}{\widetilde{\sigma}}\,\;p_x\;\;\mbox{ on }\;\; \widehat{\Sigma}_2,
\end{equation}
where $\displaystyle \{ v,p \}$ is the unique solution of (the optimality system)
\begin{equation} \label{eq3.37}
\left|
\begin{array}{l}
\displaystyle u'' - u_{xx} = 0 \ \mbox{ in } \;\; \widehat{Q}, \\
\displaystyle p'' - p_{xx} = u - \widetilde{u}_2 \ \mbox{ in } \;\; \widehat{Q},\\[5pt]
\displaystyle u = \left\{
\begin{array}{l}
\widetilde {w}_1 \ \mbox{ on } \ \widehat{\Sigma}_1,\\[5pt]
\displaystyle \frac{1}{\widetilde{\sigma}}\;\,p_x \ \mbox{ on } \ \widehat{\Sigma}_2,\\[5pt]
0 \ \mbox{ on } \ \widehat{\Sigma} \backslash \widehat{\Sigma}_0,
\end{array}
\right.\\[7pt]
p = 0 \ \mbox{ on } \ \widehat{\Sigma}, \\[5pt]
u(0) = u'(0) = 0, \\[5pt]
p(T) = p'(T) = 0, \;\; x \in \Om_t.
\end{array}
\right.
\end{equation}
Of course, $\displaystyle \{ u,p \}$ depends on $\widetilde{w}_1$:
\begin{equation}\label{cdep}
\displaystyle \{ u,p \} = \{ u(\widetilde{w}_1),p(\widetilde {w}_1) \}.
\end{equation}
\end{theorem}
\section{On the approximate controllability}\label{sec4}
\setcounter{equation}{0}
Since we have proved the existence, uniqueness and characterization of the
follower $\displaystyle \widetilde{w}_2$, the leader $\displaystyle \widetilde{w}_1$  now wants  that the solutions $u$
and $u'$, evaluated at time $t=T$, to  be as close as possible to $\displaystyle (u^0,
u^1)$. This will be possible if the system \eqref{eq3.37}  is approximately
controllable. We are looking for
\begin{equation} \label{inf1}
\begin{array}{l}
\displaystyle\inf\, \frac{1}{2\,}\,\int_{\widehat \Sigma_1} \widetilde{w}_{1}^{2}\,d{\widehat \Sigma},
\end{array}
\end{equation}
where $\displaystyle \widetilde{w}_1$ is subject to
\begin{equation} \label{subj1}
\begin{array}{l}
\displaystyle \left(u(T;{\widetilde{w}_1}), u'(T; {\widetilde{w}_1})\right)\in B_{L^2(\Om_t)}(u^0,\rho_0) \times B_{H^{-1}(\Om_t)}(u^1,\rho_1),
\end{array}
\end{equation}
assuming that  $w_1$ exists, $\rho_0$ and  $\rho_1$ being positive numbers arbitrarily
small  and $\{u^0, u^1\} \in L^2(\Om_t) \times H^{-1}(\Om_t)$.

As in \cite{Cui10}, we assume that
\begin{equation}\label{hT}
\disp T >  \frac{e^{\frac{2k(1 + k)}{(1 - k)^3}} - 1}{k}
\end{equation}
and
\begin{equation}\label{hT10}
0 < k <  1.
\end{equation}

Now as in the case \eqref{decomp 2.A} and using  Holmgren's Uniqueness
Theorem (cf. \cite{LH}; and see also \cite {Cui10} for additional discussions), the following approximate controllability result holds:

\begin{theorem}\label{AC} Assume that \eqref{hT} and \eqref{hT10} hold. Let us consider $\displaystyle \widetilde{w}_1 \in L^2(\widehat \Sigma_1)$ and $\displaystyle \widetilde{w}_2$ a Nash equilibrium in the sense \eqref{son}.
Then $\displaystyle \left(u(T), u'(T)\right)=\left(u(., T, {\widetilde{w}_1}, \widetilde{w}_2), v'(., T, {\widetilde{w}_1},
\widetilde{w}_2)\right)$, where $\displaystyle u$ solves the system \eqref{eq3.37}, generates a dense
subset of $\displaystyle L^2(\Om_t)\times H^{-1}(\Om_t)$.
\end{theorem}
\textbf{Proof.}
We decompose the solution $\displaystyle (u,p)$ of \eqref{eq3.37}
setting
\begin{equation} \label{eq3.39}
\left|
\begin{array}{l}
u = u_0 + g,\\
p = p_0 + q,
\end{array}
\right.
\end{equation}
where $u_0$, $p_0$ is given by
\begin{equation} \label{eq3.40}
\left|
\begin{array}{l}
\displaystyle u_0'' - (u_{0})_{xx} = 0 \ \mbox{ in } \ \widehat{Q},\\[5pt]\displaystyle
u_0 = \left\{
\begin{array}{l}
0 \ \mbox{ on } \ \widehat \Sigma_1,\\[5pt]
\displaystyle \frac{1}{\widetilde{\sigma}}\,({p_0})_{x} \ \mbox{ on } \ \widehat \Sigma_2,\\[5pt]\displaystyle
0 \ \mbox{ on } \ \widehat{\Sigma} \backslash \widehat\Sigma_0,
\end{array}
\right.\\[5pt]\displaystyle
u_0(0) = u_0'(0) = 0, \;\; x\in \Om_t,
\end{array}
\right.
\end{equation}
\begin{equation} \label{eq3.41}
\left|
\begin{array}{l}
\displaystyle p_0'' - (p_0)_{xx} =  u_0 - \widetilde{u}_2 \ \mbox{ in } \ \widehat Q, \\[5pt]\displaystyle
p_0 = 0 \ \mbox{ on } \ \widehat \Sigma,\\[5pt]\displaystyle
p_0(T) = p_0'(T) = 0, \;\; x \in \Om_t,
\end{array}
\right.
\end{equation}
and $\displaystyle \{g,q\}$ is given by
\begin{equation} \label{eq3.42}
\left|
\begin{array}{l}
g'' - g_{xx} = 0 \ \mbox{ in } \ \widehat Q,\\[5pt]\displaystyle
g = \left\{
\begin{array}{l}
\widetilde {w}_1 \ \mbox{ on } \ \widehat \Sigma_1,\\[5pt]
\displaystyle \frac{1}{\widetilde{\sigma}}\,q_{x} \ \mbox{ on } \ \widehat \Sigma_2,\\[5pt]\displaystyle
0 \ \mbox{on} \ \widehat \Sigma \backslash \widehat \Sigma_0,
\end{array}
\right.\\[5pt]\displaystyle
g(0) = g'(0) = 0, \;\;x \in \Om_t,
\end{array}
\right.
\end{equation}
\begin{equation} \label{eq3.43}
\left|
\begin{array}{l}
\displaystyle q'' - q_{xx} =  g \ \mbox{ in } \ \widehat Q, \\[5pt]\displaystyle
q = 0 \ \mbox{ on } \ \widehat \Sigma,\\[5pt]\displaystyle
q(T) = q'(T) = 0, \;\; x \in \Om_t.
\end{array}
\right.
\end{equation}
We next set
\begin{equation} \label{eq3.44}
\begin{array}{ccll}
A \ : & \! L^2(\widehat \Sigma_1) & \! \longrightarrow & \! H^{-1}(\Om_t) \times L^2(\Om_t) \\
& \! \widetilde{w}_1 & \! \longmapsto & \! A\,\widetilde{w}_1 = \big\{ g'(T;\widetilde{w}_1) + \delta g(T;\widetilde{w}_1),\; -g(T;\widetilde{w}_1) \big\},
\end{array}
\end{equation}
which defines
$$A \in \mathcal{L}\left( L^2(\widehat \Sigma_1); \;H^{-1}(\Om_t) \times L^2(\Om_t)\right),$$
where  $\delta$ is a positive constant.

Using \eqref{eq3.39} and \eqref{eq3.44}, we can rewrite  \eqref{subj1} as
\begin{equation} \label{subj2}
\displaystyle A\widetilde{w}_1\in \{ -u_0(T)+\delta g(T)+B_{H^{-1}(\Om_t)}(u^1,\rho_1),\;-u_0(T)+B_{L^2(\Om_t)}(u^0,\rho_0)\}.
\end{equation}
We will show that  $A\widetilde{w}_1$ generates a dense subspace of $H^{-1}(\Om_t) \times
L^2(\Om_t)$. For this, let  $\{ f^0,f^1 \} \in H_{0}^{1}(\Om_t) \times L^2(\Om_t)$ and
consider the following systems (``adjoint states"):
\begin{equation} \label{eq3.45}
\left|
\begin{array}{l}
\varphi'' - \varphi_{xx} =\displaystyle\psi \ \mbox{ in } \ \widehat Q, \\[5pt]\displaystyle
\varphi = 0 \ \mbox{ on } \ \Sigma,\\[5pt]\displaystyle
\varphi(T) = f^0, \ \varphi'(T) = f^1, \;\; x \in \Om_t,
\end{array}
\right.
\end{equation}

\begin{equation} \label{eq3.46}
\left|
\begin{array}{l}
\psi'' - \psi_{xx} = 0 \ \mbox{ in } \ \widehat Q,\\[5pt]\displaystyle
\psi = \left\{
\begin{array}{l}
0 \ \mbox{ on } \ \widehat \Sigma_1,\\[5pt]\displaystyle
\displaystyle \frac{1}{\widetilde \sigma}\,\varphi_{x} \ \mbox{ on } \ \widehat \Sigma_2,\\[5pt]\displaystyle
0 \ \mbox{ on } \ \widehat\Sigma \backslash \widehat \Sigma_0,
\end{array}
\right.\\[5pt]\displaystyle
\psi(0) = \psi'(0) = 0, \;\; x \in \Om_t.
\end{array}
\right.
\end{equation}
Multiplying \eqref{eq3.46}$_1$ by $q$, \eqref{eq3.45}$_1$ by $g$, where $q$, $g$
solve (\ref{eq3.43}) and \eqref{eq3.42}, respectively, and integrating in $Q$ we
obtain
\begin{equation} \label{eq3.47}
\int_{0}^{T}\int_{\Om_t} g\,\psi\,dx\,dt =- \frac{1}{\widetilde \sigma}\int_{\widehat \Sigma_2} q_{x}\,\varphi_{x} d\widehat \Sigma,
\end{equation}
and
\begin{equation} \label{eq3.49}
\begin{array}{l}
\displaystyle \langle g'(T),f^0 \rangle_{H^{-1}(\Om_t) \times  H_{0}^{1}(\Om_t)} + \delta \langle g(T), f^0 \rangle_{L^2(\Om_t) \times  H_{0}^1(\Om_t)}  - \big( g(T),f^1 \big)\\ = \displaystyle -\int_{\widehat \Sigma_1}\,\varphi_{x}\,\widetilde{w}_1\,d\widehat \Sigma.
\end{array}
\end{equation}
Considering the left-hand side of the above equation as the inner product of $\displaystyle \{g'(T)+
\delta g(T),-g(T)\}$ with $\{ f^0,f^1 \}$ in $ H^{-1}(\Om_t)  \times L^2(\Om_t) $ and $
H_{0}^{1}(\Om_t) \times L^2(\Om_t)$, we obtain

\begin{equation*}
\Big\langle \big\langle A\,\widetilde{w}_1 , f \big\rangle \Big\rangle = - \int_{\widehat \Sigma_1}\,\varphi_{x}\,\widetilde{w}_1\,d\widehat \Sigma,
\end{equation*}
where $\Big\langle \big\langle . , . \big\rangle \Big\rangle$ represent the duality
pairing between $ H^{-1}(\Om_t) \times L^2(\Om_t) $ and $ H_{0}^{1}(\Om_t) \times
L^2(\Om_t) $. Therefore, if
$$\langle g'(T),f^0 \rangle_{H^{-1}(\Om_t) \times H_{0}^{1}(\Om_t)} + \delta \langle g(T), f^0 \rangle_{L^2(\Om_t) \times H_{0}^1(\Om_t)} - \big( g(T),f^1 \big) = 0,$$
for all $w_1 \in L^2(\widehat \Sigma_1)$, then
\begin{equation} \label{eq3.50}
\displaystyle \varphi_{x}= 0 \ \ \mbox{ on } \ \ \widehat \Sigma_1.
\end{equation}
Hence,  in case \eqref{decomp 2.A},
\begin{equation} \label{eq3.51}
\psi = 0 \ \ \mbox{ on } \ \ \widehat \Sigma, \;\; \mbox{ so that } \psi\equiv 0.
\end{equation}
Therefore
\begin{equation} \label{eq3.54}
\begin{array}{l}
\displaystyle \varphi'' - \varphi _{xx}  = 0, \;\; \varphi = 0 \mbox{ on }  \widehat \Sigma,
\end{array}
\end{equation}
and satisfies \eqref{eq3.50}. Therefore, according to  Holmgren's Uniqueness
Theorem (cf. \cite{LH}; and see also \cite {Cui10} for additional discussions) and if \eqref{hT} holds, then $\displaystyle \varphi \equiv 0$,
so that $\displaystyle f^0=0, f^1=0$. This
ends the proof. \cqdf

\section{Optimality system for the leader}\label{sec5}
\setcounter{equation}{0}
Thanks to the results obtained in Section \ref{sec3} , we can take, for  each
$\displaystyle \widetilde{w}_1$, the Nash equilibrium $\displaystyle \widetilde{w}_2$ associated to solution $\displaystyle u$
of \eqref{eq1.3.1}. We will show the existence of a leader control $\displaystyle \widetilde{w}_1$ solution
of the following problem:
\begin{equation} \label{eq3.7cil.1}
\displaystyle \inf_{\widetilde{w}_1\in \mathcal{U}_{ad}} J(\widetilde{w}_1),
\end{equation}
where $\displaystyle \mathcal{U}_{ad}$ is the set of admissible controls
\begin{equation}\label{admcon}
\displaystyle \mathcal{U}_{ad}=\{\widetilde{w}_1\in L^2({\widehat \Sigma_1}); \; u \mbox{ solution of } \eqref{eq1.3.1} \mbox{ satisfying } \eqref{subj1}\}.
\end{equation}

For this, we will use  a duality argument due to Fenchel and Rockfellar \cite{R}
(cf. also \cite{Bre, EK}). Before, we recall the following definition from
convex analysis:
\begin{definition}\label{def123} Let E be a real linear locally convex space and  E' be its dual space. Consider any function $G: E \rightarrow ( -\infty, +\infty]$. The
function $G^{*}:  E'  \rightarrow ( -\infty, +\infty]$ defined by
$$G^{*}(x') = \sup_{x\,  \in\,  E}\{\langle x', x\rangle - G(x) \}, \;\; x' \in  E'$$
is called the conjugate function of G (sometimes called the Legendre transform
of G).
\end{definition}

The following result holds:
\begin{theorem} \label{teor3.6} Assume the hypotheses \eqref{decomp 2.A}, \eqref{hT} and \eqref{hT10}  are satisfied. Then for $\{f^0,f^1\}$
in $\displaystyle H_0^1(\Om_t) \times L^2(\Om_t)$ we uniquely define $\{\varphi, \psi, u, p \}$
by
\begin{equation} \label{eq3.139}
\left|
\begin{array}{l}
\displaystyle \varphi'' - \varphi_{xx} = \psi \ \ \text{in} \ \ \widehat Q, \\[3pt]\displaystyle
\displaystyle \psi'' - \psi_{xx}  = 0 \ \ \text{in} \ \ \widehat Q, \\[3pt]\displaystyle
\displaystyle u'' - u_{xx} = 0 \ \ \text{in} \ \ \widehat Q, \\[3pt]\displaystyle
\displaystyle p'' - p_{xx}  = u - \widetilde{u}_2 \ \ \text{in} \ \ \widehat Q, \\[3pt]\displaystyle
\displaystyle \varphi = 0 \ \ \ \text{on} \ \ \ \Sigma, \\[3pt]\displaystyle
\displaystyle \psi =
\left\{
\begin{array}{l}
\displaystyle 0 \ \ \text{on} \ \ \Sigma_1, \\[3pt]\displaystyle
\displaystyle \frac{1}{\widetilde \sigma}\;\varphi_{x} \ \ \text{on} \ \ \widehat \Sigma_2, \\[3pt]\displaystyle
\displaystyle 0 \ \ \ \text{on} \ \ \ \widehat \Sigma \backslash \widehat \Sigma_0,\\[3pt]\displaystyle
\end{array}
\right.\\
\displaystyle u =
\left\{
\begin{array}{l}
\displaystyle - \varphi_{x} \ \ \text{on} \ \ \widehat \Sigma_1,\\[3pt]\displaystyle
\displaystyle \frac{1}{\widetilde \sigma}\;p_{x} \ \ \text{on} \ \ \widehat \Sigma_2,\\[3pt]\displaystyle
\displaystyle 0 \ \ \ \text{on} \ \ \ \widehat \Sigma \backslash \widehat \Sigma_0,\\[3pt]\displaystyle
\end{array}
\right.\\
\displaystyle p = 0 \ \ \ \text{on} \ \ \ \widehat \Sigma, \\[3pt]\displaystyle
\displaystyle \varphi(.,T) = f^0,\, \varphi'(.,T) = f^1 \ \ \text{in} \ \  \Om_t, \\[3pt]\displaystyle
\displaystyle u(0) = u'(0) = 0 \ \ \text{in} \ \   \Om_t, \\[3pt]\displaystyle
\displaystyle p(T) = p'(T) = 0 \ \ \text{in} \ \  \Om_t.
\end{array}
\right.
\end{equation}
We uniquely define  $\{f^0,f^1\}$ as the solution of the variational inequality
\begin{equation} \label{eq3.140}
\begin{array}{l}
\displaystyle \big\langle u'(T,f) - u^1, \widehat{f}^0 - f^0\big\rangle_{H^{-1}(\Om_t) \times  H_{0}^{1}(\Om_t)} - \big(u(T,f) - u^0, \widehat{f}^1 - f^1\big)  \\[10pt]
\displaystyle + \rho_1\big(||\widehat{f}^0|| - ||f^0||\big) + \rho_0\big(|\widehat{f}^1| - |f^1|\big) \geq 0,\,\forall\, \widehat{f} \in H_{0}^{1}(\Om_t) \times L^2(\Om_t).
\end{array}
\end{equation}
Then the optimal leader is given by
\begin{equation*}
\widetilde{w}_1 = - \varphi_{x} \ \ \text{on} \ \  \widehat \Sigma_1,
\end{equation*}
where $\varphi$ corresponds to the solution of \eqref{eq3.139}.
\end{theorem}
\textbf{Proof.} We introduce two convex proper functions as follows, firstly
\begin{equation}\label{eq3.119}
\begin{array}{l}
\displaystyle F_1 : L^2(\widehat \Sigma_1) \longrightarrow \re \cup \{\infty\},\\[5pt]
\displaystyle F_1(\widetilde{w}_1) = \frac{1}{2} \int_{\widehat \Sigma_1}\widetilde{w}_{1}^{2}\,d\widehat \Sigma
\end{array}
\end{equation}
the second one
\begin{equation*}
F_2 : H^{-1}(\Omega_t) \times L^2(\Omega_t) \longrightarrow \re \cup \{\infty\},
\end{equation*}
given by
\begin{equation} \label{eq3.120}
\begin{array}{ccl}
F_2(A\widetilde{w}_1) = F_2\big(\{g'(T,\widetilde{w}_1) + \delta g(T,\widetilde{w}_1),-g(T,\widetilde{w}_1)\}\big) \\
= \left\{
\begin{array}{l}
0, \text{ if }
\left\{
\begin{array}{l}
g'(T) + \delta g(T) \in u^1 - u_0'(T) + \delta g(T) + \rho_1B_{H^{-1}(\Om_t)},\\
-g(T) \in -u^0 + u_0(T,\widetilde{w}_1) - \rho_0B_{L^2(\Om_t)},
\end{array}
\right.\\
+ \infty, \text{ otherwise}.
\end{array}
\right.
\end{array}
\end{equation}
With these notations, problems \eqref{inf1}--\eqref{subj1} become equivalent to
\begin{equation} \label{eq3.122}
\begin{array}{l}
\displaystyle \inf_{\widetilde{w}_1 \in L^2(\widehat \Sigma_1)}\big[F_1(\widetilde{w}_1) + F_2(A\widetilde{w}_1)\big]
\end{array}
\end{equation}
provided we prove that the range of $\displaystyle A$ is dense in $\displaystyle H^{-1}(\Om_t)
\times L^2(\Om_t)$, under conditions \eqref{hT} and  \eqref{hT10}.

By the Duality Theorem of Fenchel and Rockafellar \cite{R}(see also \cite{Bre, EK}), we have
\begin{equation} \label{eq3.124}
\begin{array}{l}
\inf_{\widetilde{w}_1 \in L^2(\widehat \Sigma_1)}[F_1(\widetilde{w}_1) + F_2(A\widetilde{w}_1)]\\[5pt]\displaystyle  = -\inf_{(\widehat{f}^0,\widehat{f}^1) \in H_{0}^{1}(\Om_t) \times L^2(\Om_t)} [F_{1}^{*}\big(A^*\{\widehat{f}^0,\widehat{f}^1\}\big) + F_{2}^{*}\{-\widehat{f}^0, -\widehat{f}^1\}],
\end{array}
\end{equation}
where $\displaystyle F_i^*$ is the conjugate function of $\displaystyle F_i  (i=1,2)$ (see Definition \ref{def123})
and $\displaystyle A^*$ the adjoint of $\displaystyle A$.

We have
\begin{equation} \label{eq3.121}
\begin{array}{ccccc}
A^* \ : & \! H_{0}^{1}(\Omega_t) \times L^2(\Omega_t) & \! \longrightarrow & \! L^2(\widehat \Sigma_1) \\
& \! (f^0,f^1) & \! \longmapsto & \! A^*f = & \! - \varphi_{x},
\end{array}
\end{equation}
where $\varphi$ is given in (\ref{eq3.45}).
We see easily that
\begin{equation} \label{eq3.125}
F_{1}^{*}(\widetilde{w}_1) = F_1(\widetilde{w}_1)
\end{equation}
and
\begin{equation}\label{eq3.125.2}
\begin{array}{ccl}
\displaystyle  F_{2}^{*}(\{\widehat{f}^0,\widehat{f}^1 \})  &=& \displaystyle \langle u^1 - u_0'(T) + \delta g(T), \widehat{f}^0\rangle_{H^{-1}(\Omega_t) \times H_{0}^{1}(\Omega_t)} \\&&+\big( u_0(T) - u^0 ,\widehat{f}^1\big)
\displaystyle + \rho_1||\widehat{f}^0|| + \rho_0|\widehat{f}^1|.
\end{array}
\end{equation}
Therefore the (opposite of) right-hand side of  \eqref{eq3.124} is given by
\begin{align} \label{eq3.127}
\displaystyle
& - \inf_{\widehat{f} \in H_{0}^{1}(\Omega_t) \times L^2(\Omega_t)} \bigg\{\frac{1}{2}\int_{\widehat \Sigma_1}{\varphi}_{x}^2d\widehat \Sigma  + \big( u^0 - u_0(T) ,\widehat{f}^1\big) \\
\nonumber & - \langle u^1 - u_0'(T) + \delta g(T), \widehat{f}^0\rangle_{H^{-1}(\Omega_t) \times H_{0}^{1}(\Omega_t)} + \rho_1||\widehat{f}^0|| + \rho_0|\widehat{f}^1|\bigg\}.
\end{align}
This is the dual problem of  \eqref{inf1},  \eqref{subj1}.

We  have now two ways to derive the optimality system for the leader control,
starting from the primal or from the dual problem. \cqdf

\begin{remark} The main results of this paper could be extended to the case of the exact controllability. The key point is to define directly the energy function of a wave equation in the non-cylindrical domain and use the multiplier method to overcome difficulties. For interested readers on this subject, we cite for instance \cite{Cui10} .
\end{remark}


\paragraph{\bf Acknowledgements}

The author wants to express his gratitude to the
anonymous reviewers for their questions and commentaries; they were very
helpful in improving this article.

\paragraph
\end{document}